\documentclass[12pt,leqno]{article}
\usepackage{amssymb}
 \usepackage[french]{babel}
\newcommand{\dd}{{\rm \kern 3pt I\kern-9pt d}}

\newcommand{\Abar}{{\backslash\kern-8pt A}}

\newtheorem{Th}{Theorem}

\newtheorem{Pro}[Th]{Proposition}

\newcommand{\CC}{\mbox{{\helv \i}{\rm\kern-5pt C}}}

\begin{document}
\begin{center}\bf\large
UNE STRUCTURE UNIFORME SUR UN ESPACE $\mathcal{F}(E,F)$\\
\end{center}
\begin{center}
Cahiers de Topologie et G\'eom\'etrie diff\'erentielle

Vol XI, 2, (1969) p207-214
\end{center}
\begin{flushright}
Nicolas BOULEAU
\end{flushright}

Soient $E$ un espace topologique, et $F$ un espace uniforme. L'objet de cette \'etude est d'introduire une topologie sur $\mathcal{F}(E,F)$, ensemble des applications de $E$ dans $F$, telle que toute limite d'une suite de fonctions continues soit continue et que, si une suite de fonctions continues converge en chaque point vers une fonction continue, elle converge vers cette fonction pour cette topologie. 

Plus pr\'ecis\'ement, nous d\'efinirons sur $\mathcal{F}(E,F)$ une structure uniforme poss\'edant les propri\'et\'es suivantes: 

1) $\mathcal{C}(E,F)$ est ferm\'e dans $\mathcal{F}(E,F)$ muni de cette structure uniforme. 

2) La restriction de cette structure uniforme \`a $\mathcal{C}(E,F)$ est \'equivalente \`a la structure uniforme de la convergence simple. 
Nous \'etudions ensuite plusieurs applications de cette topologie. 

\section{Structure uniforme de la $V$-convergence}
Soit $E$ un espace topologique et soit $F$ un espace uniforme. Soit $W$ un entourage de $F$, $A$ une partie finie de $E$ et soit $\mathcal{U}_{W,A}$ l'ensemble des couples $(f,g)\in\mathcal{F}(E,F)\times\mathcal{F}(E,F)$  tels qu'il existe des voisinages $V_{a_1},\ldots,V_{a_n}$ des points de $A$ tels que l'on ait :
$$\forall x\in\cup_{i=1}^n V_{a_i}\quad\quad(f(x),g(x))\in W.$$
\begin{Pro}{Lorsque $W$ d\'ecrit un syst\`eme fondamental d'entourages de $F$ et lorsque $A$ d\'ecrit l'ensemble des parties finies de $E$, les $\mathcal{U}_{W,A}$ d\'ecrivent un syst\`eme fondamental d'entourages d'une structure uniforme sur $\mathcal{F}(E,F)$, appel\'ee structure uniforme de la $V$-convergence. }
\end{Pro}
DEMONSTRATION. il suffit de v\'erifier que les axiomes des syst\`emes fondamentaux d'entourages sont v\'erifi\'es (cf. \cite{bourbaki1} Chap. 2 Structures Uniformes). Or : 
$$W_3\subset W_1\cap W_2\Rightarrow\mathcal{U}_{W_3,A_1\cup A_2}\subset\mathcal{U}_{W_1,A_1}\cap\mathcal{U}_{W_2,A_2},$$
$$W^\prime\subset \stackrel{-1}{W}\Rightarrow \mathcal{U}_{W^\prime,A}\subset\stackrel{-1}{\mathcal{U}}_{W,A}$$
$$\stackrel{2}{W}_1\subset W\Rightarrow\stackrel{2}{ \mathcal{U}}_{W_1,A}\subset\mathcal{U}_{W,A},$$
et, comme chaque $\mathcal{U}_{W,A}$ contient la diagonale de $\mathcal{F}(E,F)$, la proposition est d\'emontr\'ee. \\

\noindent$\mathcal{F}(E,F)$ muni de cette structure sera not\'e $\mathcal{F}_V(E,F)$. 

\begin{Pro} La structure uniforme de la $V$-conuergence est plus fine que celle de la convergence simple et moins fine que celle de la conuergence uniforme locale. 
\end{Pro}
Ces propri\'et\'es se voient imm\'ediatement en comparant les filtres 
d'entourages de ces diff\'erentes structures uniformes. 
\begin{Pro} L'ensemble $\mathcal{C}(E,F)$ est ferm\'e dans $\mathcal{F}(E,F)$. Les structures uniformes induites sur $\mathcal{C}(E,F)$ par les structures uniformes de la convergence simple et de la $V$-convergence sont \'equivalentes. 
\end{Pro}
DEMONSTRATION. 

a) Soit $g\in\overline{\mathcal{C}(E,F)}^V$ adh\'erence de $\mathcal{C}(E,F)$ dans $\mathcal{F}_V(E,F)$, soit $x_0\in E$,  montrons que $g$ est continue en $x_0$ :

 Soient $W$ et $W^\prime$ des 
entourages de $F$ tels que $\stackrel{3}{W'} \subset W$. Il existe $f\in\mathcal{C}( E , F)$ telle que $(f. g)\in\mathcal{U}_{W^\prime, x_0}$, c'est-\`a-dire qu'il existe un voisinage $V^1_{x_0}$ de $x_0$ tel que: 
	$$\forall x\in V^1_{x_0}, \qquad(f(x),g(x))\in W^\prime.$$
D'autre part $f$  \'etant continue, il existe un voisinage $V^2_{x_0}$ de $x_0$  tel que: 
$$\forall x\in V^2_{x_0}, \qquad(f(x),f(x_0))\in W^\prime,$$ donc
$$\forall x\in V^1_{x_0}\cap V^2_{x_0}, \qquad(g(x),g(x_0))\in \stackrel{3}{W^\prime}\subset W.$$
donc  $g$ est continue en $x_0.$

b) Il suffit de montrer que sur $\mathcal{C}( E, F)$ la structure uniforme de la $V$-convergence est moins fine que celle de la convergence, simple: i.e. pour tout $\mathcal{U}_{W,A}$ entourage de $\mathcal{C}_V(E,F)$ il existe un entourage $\mathcal{W}$ de $\mathcal{C}_s(E,F)$ tel que $\mathcal{W}\subset\mathcal{U}_{W,A}$. Soit donc
$$\mathcal{U}_{W,A}=\{(f,g) : \exists V_{a_1},\ldots,V_{a_n} : \forall x\in\cup_{i=1}^n V_{a_i},(f(x),g(x))\in \mathcal{W}\}$$
et soit $W'$ un entourage de $F$ tel que $\stackrel{3}{W'}\subset W$. Consid\'erons l'entourage de $\mathcal{C}_s(E,F)$ d\'efini par
$$\mathcal{W}_{W',A}=\{(f,g) : \forall a_i\in A, (f(a_i),g(a_i))\in W'\}.$$ Le fait que $f$ et $g$ soient continues implique que $\mathcal{W}_{W',A}\subset\mathcal{U}_{W,A}$.\\

\noindent REMARQUES

1)  Si la structure uniforme de $F$ est d\'efinie par les serni-distances $d_i$, $i\in I$, la structure uniforme de la $V$áconvergence peut se d\'efinir par les semi-distances : 
$$\delta_{i,A}(f,g)=\inf_{V_k\in\mathcal{V}_{a_k}}\sup_{x\in\cup_kV_k}d_i(f(x),g(x))$$
o\`u $\mathcal{V}_{a_k}$ d\'esigne l'ensemble des voisinages de $a_k$ dans $E$ et $A$ l'ensemble fini $\{a_1,\ldots,a_n\}$.

2) Si un filtre $\mathfrak{F}$ converge dans $\mathcal{F}_V(E.F)$ : $f=\lim_{\mathfrak{F}} f_i$ alors 
pour tout $x$ dans $E$ et tout entourage $W$ de $F$, il existe un 	$V\in \mathcal{V}_x$
(ensemble des voisinages de $x$) et une fonction $f_i\in\mathfrak{F}$ tels que $f$  et $f_i$ soient voisins d'ordre $W$ dans $V$. 
\begin{Pro} Soient $E$ un espace topologique et $F$ un espace uniforme s\'epar\'e, soit $H\subset\mathcal{C}( E. F)$ pour que $H$ soit relativement compact pour la $V$-convergence, il faut et il suffit que : 

1) $\forall x$  $H( x)$ soit relativement compact dans $F$ et 

2) $\overline{H}^s\subset\mathcal{C}( E. F)$, $\overline{H}^s$ d\'esignant la fermeture de $H$ pour la topologie de la convergence simple. 
\end{Pro}
DEMONSTRATION. 

D'apr\`es la proposition 3, ces conditions sont suffisantes. Montrons qu'elles sont nŽcessaires: 

Supposons donc $H$ relativement compact pour la $V$-convergence, $\forall x\in E$ l'application $f \mapsto f( x)$ de $\mathcal{F}_V( E. F)$ dans $F$ est continue, puisque la $V$-convergence est plus fine que la convergence simple, donc $H(x)$ est relativement compact dans $F$. 

Pour la m\^eme raison $\overline{H}^V\subset\overline{H}^s$ et, puisque l'injection canonique de $\mathcal{F}_V(E ,F)$ dans $\mathcal{F}_s(E ,F)$ est continue, $\overline{H}^V$ est compact dans $\mathcal{F}_s(E ,F )$, donc ferm\'e. Comme $\overline{H}^V$ contient $H$, n\'ecessairement $\overline{H}^V$ contient $\overline{H}^s$, donc $\overline{H}^V=\overline{H}^s$. D'apr\`es la proposition 3, on a donc: $\overline{H}^s\subset\mathcal{C}(E,F).$
\section{Crit\`ere de $V$-convergence}
M\^eme si $F$ est complet,  $\mathcal{F}_V(E ,F)$ n'est pas n\'ecessairement complet; toutefois il existe un crit\`ere de convergence qui, comme celui de Cauchy, ne fait pas intervenir la limite du filtre \'etudi\'e: 
\begin{Pro} Soit $E$ un espace topologique et soit $F$ un espace uniforme complet d\'efini par les semi-distances $d_i,\;i\in I$. Une condition n\'ecessaire et suffisante pour qu'un filtre $\mathfrak{F}$ sur $\mathcal{F}_V(E ,F)$ soit convergent est que: 
$$\forall i\in I, \forall\varepsilon>0, \forall a\in E, $$
$$\exists A\in\mathfrak{F} : \forall f\in A, \exists V_a : \forall x\in V_a, \exists B\in\mathfrak{F} :\forall g\in B$$
$$d_i(f(x),g(x))\leq \varepsilon.$$
\end{Pro}
DEMONSTRATION

{\it Condition n\'ecessaire.} Si $\mathfrak{F}$ $V$-converge vers $\ell$, on a:
$$\forall i\in I, \forall\varepsilon>0, \forall a\in E, $$
$$\exists A\in\mathfrak{F} : \forall f\in A, \exists V_a : \forall x\in V_a : d_i(f(x),\ell(x))\leq\frac{\varepsilon}{2}.$$
Soit $x$ fix\'e, le filtre $\mathfrak{F}(x)$ converge vers $\ell(x)$ dans $F$, donc $\exists B : \forall g\in B, d_i(g(x),\ell(x))\leq \frac{\varepsilon}{2}$, d'o\`u le r\'esultat.

{\it Condition suffisante.} On remarque d'abord qu'elle implique que le filtre $\mathfrak{F}( x)$ est un filtre de Cauchy dans $F$, soit $\ell(x)$ sa limite, alors on voit qu'elle implique que le filtre $\mathfrak{F}$  $V$-converge vers $\ell$.\\
 
\noindent REMARQUE. Si $F$ est seulement s\'equentiellement complet, on a un crit\`ere analogue pour les suites. Par exemple, si $Y$ est un espace vectoriel norm\'e s\'equentiellement complet, une suite $\{f_n\}$ converge dans $\mathcal{F}_V(E ,F)$  si, et seulement si:
$$\forall\varepsilon>0,\forall a\in X, \exists N :\forall n\geq N, \exists V^n_a :\forall x\in V^n_a, \exists P: \forall p\geq P : \|f_n(x)-f_p(x)\|\leq\epsilon.$$
On d\'eduit de ce crit\`ere les applications suivantes: 
\begin{Pro}
 Soient $X$ un espace topologique, $Y$ un espace de Banach. Soit $\{f_n\}$ une suite d'applications continues de $X$ dans $Y$. On suppose 
que la s\'erie $\sum_n\|f_n\|$ converge en chaque point $x$ vers une fonction 
$\Sigma(x)$ continue. Alors la s\'erie $\sum_nf_n(x)$ converge vers $S ( x)$ continue. 
\end{Pro}
DEMONSTRATION. Comme $\|\sum_{k=p}^{k=q}f_k(x)\|\leq \sum_{k=p}^{k=q}\|f_k(x)\|$, le fait que 
	la sŽrie $\sum\|f_n(x)\|$  $V$-converge entra\^ine, d'apr\`es le crit\`ere ci-dessus, que 
la s\'erie $\sum_nf_n(x)$ $V$-converge. 
\begin{Pro}Soit  $\{f_n\}$  une suite d'applications continues de $X$ dans un espace 
de Banach $Y$ v\'erifiant $\|\sum_{k=0}^nf_k(x)\|\leq A(x)$, $ A ( x)$ \'etant localement 
born\'ee. Soient $\{\varepsilon_n(x)\}$ continues de $X$ dans $\mathbb{R}$ tendant vers z\'ero pour $n\rightarrow\infty$, 
 et telles que la s\'erie $\sum_n|\varepsilon_n(x)-\varepsilon_{n-1}(x)|$  converge uers une fonction 
continue. Alors la s\'erie $\sum_nf_n(x)\varepsilon(x)$ converge vers une fonction continue de $X$ dans $Y$. \end{Pro}
Cette proposition se d\'emontre comme on d\'emontre habituellement 
la r\`egle d'Abel pour les s\'eries en utilisant le crit\`ere ci-dessus au lieu du 
crit\`ere de Cauchy. 
\section{Sous-ensembles ferm\'es pour la $V$-convergence et propri\'et\'es uniformes semi-locales}
Soient toujours $E$ un espace topologique et $F$ un espace uniforme; on dira qu'une propri\'et\'e $P$ ($P \subset\mathcal{F} (E, F))$ est uniforme semi-locale si, et seulement si, pour tout $f\in\mathcal{F} (E, F)$ la condition : 
$$\forall x\in E, \forall W \mbox{ entourage de }F, \exists g\in P \mbox{ et } \exists V\in \mathcal{V}_x \mbox{ avec } f \mbox{ et } g \mbox{ voisines d'ordre } W \mbox{ sur } V$$
implique $f\in P$.
\begin{Pro} Pour que  $P \subset\mathcal{F} (E. F))$ soit une propri\'et\'e uniforme semi-locale, il faut et il suffit que $P$ soit ferm\'e pour la $V$-convergence. 
\end{Pro}
La condition ci-dessus est en effet \'equivalente \`a la suivante: 
$$\forall A\subset  E, A \mbox{ fini },\forall W \mbox{ entourage de }F, \exists g\in P \mbox{ et } \exists V_n\in \mathcal{V}_{a_n} \mbox{ avec } f \mbox{ et } g \mbox{ d'ordre } W \mbox{ sur } \cup_nV_n$$
qui peut elle-m\^eme s'\'enoncer: 
$$\forall \mathcal{W}_{W,A} \mbox{ entourage de }\mathcal{F}_V(E,F), \exists g\in P \mbox{ tel que } f \mbox{ et } g \mbox{ soient d'ordre } \mathcal{W}_{W,A} $$
c'est-\`a-dire $f\in \overline{P}$ d'o\`u la proposition. \\

Dans la suite nous allons rechercher quelques propri\'et\'es uniformes semi-locales. Exemples : continuit\'e, si $F = \mathbb{R}$ semi-continuit\'e inf\'erieure ou sup\'erieure. De m\^eme on montre facilement le r\'esultat suivant: 
\begin{Pro}Soient $X$ un espace topologique et $Y$ un espace uniforme m\'etrisable, soit $a_n\in X$  avec $\lim a_n=a$. Soit $H_{\{a_n\}}$ l'ensemble des fonctions de $X$ dans $Y$ telles que $\lim_n f( a_n )$ existe. Alors: 

a) $H_{\{a_n\}}$ est ferm\'e dans $\mathcal{F}_V(X. Y)$; 

b) si une suite de fonctions $f_p\in H_{\{a_n\}}$ $V$-conuerge vers $g$ on a : 
$$\lim_p(\lim_n f_p(a_n))=\lim_ng(a_n).$$
	\end{Pro}
En particulier l'ensemble des fonctions r\'egl\'ees de $\mathbb{R}$ dans $\mathbb{R}$ est ferm\'e pour la $V$-convergence. 
\begin{Pro} Soit $X$ un espace topologique et soit $Y$ un e.v.t. Le sous-espace vectoriel $\mathcal{F}_{lb}(X,Y)$ de  $\mathcal{F}(X,Y)$ des fonctions localement born\'ees est ferm\'e dans $\mathcal{F}_{V}(X,Y)$.
La structure de la $V$-convergence est compatible avec la structure d'espace vectoriel de ce sous-espace et est localement convexe si $Y$ l'est. \end{Pro}
DEMONSTRATION 

a) Etre localement born\'e est une propri\'et\'e uniforme semi-locale donc $\mathcal{F}_{lb}(X,Y)$ est ferm\'e.

b) Soit $W$ un voisinage de $0$ dans $Y$; les ensembles 
$$\mathcal{W}_{W,A} =\{f : \exists V_n\in\mathcal{V}_{a_n} : \forall y\in\cup_nV_n, f(y)\in W\}$$
forment, lorsque $W$ d\'ecrit l'ensemble des voisinages \'equilibr\'es de $0$ dans $Y$, un syst\`eme invariant par homoth\'etie et form\'e de voisinages \'equilibr\'es. 

Comme on voit facilement que la $V$-convergence est compatible avec la structure de groupe additif de $\mathcal{F}_{lb}(X, Y)$ et comme les $f\in\mathcal{F}_{lb}(X, Y)$ sont localement born\'ees et donc les $\mathcal{W}_{W,A} \cap \mathcal{F}_{lb}(X, Y)$  absorbants, la structure de la $V$-convergence est compatible avec la structure d'espace vectoriel de $\mathcal{F}_{lb}(X, Y)$. La derni\`ere assertion r\'esulte de ce que, si $W$ est convexe, $\mathcal{W}_{W,A}$ l'est \'egalement. 
\begin{Pro} Soit $X$ localement compact et soit $\mu$ une mesure de Radon positive sur $X$. L'ensemble des fonctions $\mu$-mesurables et l'ensemble des fonctions localement $\mu$-int\'egrable  sont ferm\'es pour la $V$-convergence. \end{Pro}
DEMONSTRATION

Cela revient \`a montrer que les propri\'et\'es d'\^etre $\mu$-mesurable et d'\^etre localement $\mu$-int\'egrable sont des propri\'et\'es uniformes semi-Iocales. Soit donc $f\in\mathcal{F}( X, \mathbb{R})$ telle que : 
$$\forall x\in X, \forall\varepsilon>0, \exists V\in \mathcal{V}_x \mbox{ et }\exists f_{\varepsilon, V} \mbox{ telle que }|f-f_{\varepsilon,V}|\leq \varepsilon \mbox{ dans }V,$$
$$f_{\varepsilon,V}\mbox{ \'etant }\mu-\mbox{mesurable (resp. localement }\mu\mbox{-int\'egrable)}.$$
Soit $K$ un compact de $X$. Il existe un nombre fini de $V$ soient $V_1,\ldots,V_n$ recouvrant $K$. La fonction $g_\varepsilon$ d\'efinie par 
$$g_\varepsilon=f_{\varepsilon,V_k}\mbox{ sur }V_k\backslash(\cup_{i=1}^{k-1}V_i)$$
est $\mu$-mesurable (resp. localement $\mu$-int\'egrable) si on choisit les $V_k$ ouverts, et l'on a $|g_\varepsilon-f|\leq\varepsilon$ dans $K$. Donc $f$ est $\mu$-mesurable (resp. localement $\mu$-int\'egrable).\\

Si $Y$ est m\'etrisable, $\mathcal{F}_V( X, Y)$ n'est pas n\'ecessairement m\'etrisable; toutefois on a : 
\begin{Pro} Soit $X$ un espace localement compact d\'enombrable \`a l'infini et soit $Y$ un espace uniforme m\'etrisable. Soit $A\subset\mathcal{F}_V( X, Y)$ ; si $f_0\in\overline{A}$,  il existe un sous-ensemble d\'enombrable $A_1$ de $A$ tel que $f_0\in\overline{A_1}$.
\end{Pro}
DEMONSTRATION

Soit $d$ une distance sur $Y$ d\'efinissant la structure de $Y$. Soit $K$ un compact de $X$; on voit en utilisant la compacit\'e de l'espace $K^m$ que, pour tout couple d'entiers $m , n$ , il existe une partie finie $A_{m,n}$ 
de $A$ telle que, pour tout ensemble de $m$ points $t_k$ de $K$, il existe $f\in A_{m,n}$ et des voisinages $V_k\in\mathcal{V}_{t_k}$ tels que $d(f_0,f)\leq\frac{1}{n}$ sur $\cup_{k=1}^mV_k$.  D'o\`u on 
d\'eduit la propri\'et\'e, puisque $X$ est r\'eunion d\'enombrable de compacts. 
\begin{Pro} Soient $X,Y$ deux espaces m\'etrisables, $X$ \'etant localement compact d\'enombrable \`a l'infini. Alors

a) l'ensemble des fonctions bor\'eliennes est ferm\'e dans $\mathcal{F}_V(X, Y)$. 

b) Pour tout ordinal d\'enombrable $\alpha$, l'ensemble des fonctions bor\'eliennes de classe $\alpha$ est ferm\'e dans $\mathcal{F}_V(X, Y)$. 
\end{Pro}
DEMONSTRATION

Cette proposition se d\'emontre exactement comme la proposition 11: Avec les m\^emes notations $g_\varepsilon$ est bor\'elienne (resp. bor\'elienne de classe $\alpha$) d\`es que les $f_{\varepsilon,V_k}$ sont bor\'eliennes (resp. bor\'eliennes 
de classe $\alpha$) et les $V_k$ ferm\'es. Comme $g_\varepsilon$ approche $f$ uniform\'ement sur 
$K$, on en d\'eduit que la restriction de $f$ \`a $K$ est bor\'elienne (resp. bor\'elienne de classe $\alpha$). Il en r\'esulte, puisque $X$ est r\'eunion d\'enombrable de compacts donc de ferm\'es, que $g$ est bor\'elienne (resp. bor\'elienne de classe $\alpha$). \\

Les deux propositions suivantesá sont les expressions en termes de $V$-convergence de propri\'et\'es de fonctions continues. 
\begin{Pro} Soit $H\subset\mathcal{F}(X, Y)$, $X$ compact, $Y$ m\'etrisable, $H$ \'etant form\'e d'appli\-cations continues et tel que toute suite de points de $H$ admette une valeur d'adh\'erence pour la $V$-convergence. Alors $\overline{H}^V$ est compact. \end{Pro}
DEMONSTRATION Cf. Bourbaki \cite{bourbaki2} Chap. IV,\S 2, ex. 15.\\

\noindent{\bf Th\'eor\`eme de Stone-Weierstrass}. {\it Soit $X$ un espace compact et $H\subset\mathcal{C}(X,\mathbb{R})$ tel que :
$$(u\in H\mbox{ et } v\in H)\Rightarrow(\sup(u,v)\in H\mbox{ et }\inf(u,v)\in H).$$
Alors $\overline{H}^u=\overline{H}^V$ o\`u $\overline{H}^u$ d\'esignant la fermeture de $H$ pour la convergence uniforme.}\\

\vspace{1cm}

Nicolas BOULEAU

51 rue G\'erard

75013 Paris

\end{document}